\input amstex

\define\pf{\demo{Proof}}
\define\endpf{\enddemo}
\define\endst{\endproclaim}
\define\thm{\proclaim{Theorem}}
\define\Thm#1{\proclaim{#1.  Theorem}}
\define\thM#1{\proclaim{Theorem \rm(#1)}}
\define\lemma{\proclaim{Lemma}}
\define\prop{\proclaim{Proposition}}
\define\Prop#1{\proclaim{#1.  Proposition}}

\define\rem{\remark{Remark}}
\define\rems{\remark{Remarks}}
\define\endrem{\endremark}

\define\endnote{\endremark}
\define\defn{\definition{Definition}}
\define\Defn#1{\definition{#1.  Definition}}
\define\enddefn{\enddefinition}

\documentstyle{amsppt}

\hcorrection{.25in}

\document

\magnification=1200
\baselineskip=15pt

\font\ff=cmbx10 at 14 truept
 at 16 truept
\font\ffl=cmbx10 
\font\ffi=cmbxti10

\def\leaderfill{\leaders\hbox to .35em{\hss.\hss}\hfill}

\redefine\d{\Cal D}
\redefine\O{\Cal O}
\define\M{\Cal M}
\redefine\H{\Cal H}
\redefine\I{\Cal I}
\redefine\T{{\Cal T}_\ell}

\define\Q{{\Bbb Q}}
\define\Qbar{{\overline {\Bbb Q}}}
\define\R{{\Bbb R}}
\define\C{{\Bbb C}}
\define\Z{{\Bbb Z}}
\define\N{{\Bbb N}}
\define\k{{\Bbbk}}
\define\Gm{{{\Bbb G}_m}}

\define\Zl{{\Z_\ell}}
\define\Ql{{\Q_\ell}}
\define\Cl{{\C_\ell}}
\define\Qlbar{{\overline \Q_\ell}}

\redefine\sp{\frak {sp}}
\redefine\sl{\frak {sl}}
\redefine\gl{\frak {gl}}
\redefine\u{\frak u}

\redefine\so{{\frak s}{\frak o}}
\define\h{\frak h}
\define\hh{\widetilde \h}
\define\g{\frak g}
\redefine\gg{\widetilde \g}
\define\a{\frak a}
\redefine\b{\frak b}
\redefine\i{\frak i}
\redefine\l{\frak l}
\redefine\ll{\widetilde \l}

\define\G{\frak G}

\define\p{\wp}

\redefine\embed{\hookrightarrow}
\define\iso{\cong}
\redefine\imply{\Rightarrow}

\redefine\wedge{\bigwedge}

\define\av{abelian variety}
\define\avs{abelian varieties}
\define\rep{re\-pre\-sen\-ta\-tion}
\define\reps{re\-pre\-sen\-ta\-tions}
\define\irr{irreducible}
\define\irrep{\irr\ \rep}
\define\irreps{\irr\ \reps}
\redefine\ss{semi-simple}
\define\et{{\text{\'et}}}
\define\wrt{with res\-pect to}

\define\QED{\hfill$\square$}

\define\HTd{the Hodge-Tate decomposition}
\define\HT{Hodge-Tate}
\define\Hd{Hodge decomposition}
\define\ext{extension}
\define\red{reduction}
\define\s.t.{such that}
\define\nf{number field}

\define\loc.cit.{{\it loc.\ cit.}}

\def\negspi{\kern-.06em}
\def\negspii{\kern-.15em}

\define\Lie{\operatorname{{\Cal L}\negspi\text{\it ie}}}

\define\Gal{\operatorname{{\Cal G}\!\text{\it a}\ell}}

\define\Hg{\operatorname{{\Cal H}\negspii\text{\it g}}}

\define\End{\operatorname{\text{\rm End}}}
\define\GL{\operatorname{\text{\rm GL}}}
\define\Aut{\operatorname{\text{\rm Aut}}}
\define\Ker{\operatorname{\text{\rm Ker}}}
\define\rk{\operatorname{\text{\rm rk}}}


\catcode`\@=11
\def\logo@{\baselineskip2pc \hbox to\hsize{\hfil\eightpoint }}
\catcode`\@=\active


\topmatter

\title
{On the Mumford--Tate Conjecture \\
for Abelian Varieties \\
with Reduction Conditions}
\endtitle

\author
{\bf Alex Lesin}
\endauthor

\thanks
{Research at MSRI supported in part by N.S.F. grant \#DMS 9022140.}
\endthanks
\address
Alex Lesin,
Mathematical Sciences Research Institute,
1000 Centennial Drive,
Berkeley, CA 94720
\endaddress
\email
shura\@msri.org
\endemail

\leftheadtext\nofrills{A. Lesin}
\rightheadtext\nofrills{On the Mumford--Tate Conjecture}



\abstract
We study monodromy action on abelian varieties satisfying certain bad reduction
conditions. These conditions allow us to get some control over the Galois
image.
As a consequence we verify the Mumford--Tate conjecture for such abelian
varieties.
\endabstract

\toc\nofrills{{\ffl Contents}}

\widestnumber\head{\bf II}

\widestnumber\subhead{10}

\specialhead{}
{Introduction}
\endspecialhead

\head 
{\bf I.} {\bf Preliminaries, first applications}
\endhead

\subhead 
0.  Basics and notations
\endsubhead
\subhead 
1.  On the Hodge-Tate decomposition
\endsubhead
\subhead 
2.  Abelian 4-folds 
\endsubhead

\head
{\bf II.} {\bf Abelian varieties with reduction conditions}
\endhead

\subhead
3.  Bad reduction and monodromy action
\endsubhead
\subhead
4.  Minimal reduction
\endsubhead
\subhead
5.  Applications of minimal reduction
\endsubhead
\subhead
6.  Another type of bad reduction
\endsubhead
\subhead
7.  A curious result
\endsubhead

\specialhead{}
{References}
\endspecialhead

\endtoc

\endtopmatter

\head
{\ff Introduction}
\endhead

There are two long outstanding conjectures
due to Hodge and Tate related to the structure of the ring of algebraic
cycles modulo homological equivalence. The Mumford--Tate conjecture
implies that for abelian varieties the two are equivalent.
The focus of this work is on the Mumford--Tate conjecture for special classes of
abelian varieties.

Let $X$ be a smooth projective algebraic variety over $\C$.
Then its $r$-th cohomology group admits the Hodge decomposition:
$H^r(X, \C) = \underset{p+q=r} \to \oplus H^{p,q}(X)$. The {\it Hodge cycles}
are those rational cohomology classes, i.e., elements of $H\,\dot{ }(X, \Q)$
that sit in the components $H^{p,p}(X)$ via the canonical embedding $H\,\dot{
}(X, \Q) \embed H\,\dot{ }(X, \C)$. Denote $\H^p(X) := H^{p,p}(X) \cap
H^{2p}(X, \Q)$ the group of codimension $p$ Hodge cycles.
Then that $\H(X) := \oplus_p \H^p(X)$ has a ring structure with respect to the
cup-product. It is immediate that rational linear combinations of the
cohomology classes of algebraic subvarieties in $X$ (=: {\it algebraic cycles})
are Hodge.
The {\it Hodge conjecture} claims that the converse is also true, viz., {\it
all the Hodge cycles are algebraic}.

The only general result in this direction is the Lefschetz (1,1)-theorem
asserting algebraicity of all codimension 1 Hodge cycles (= rational (1,1)
cohomology classes, hence the name).
Denote $\d(X)$ the subring of $\H(X)$ generated by $\H^1,\ \d^p(X) := \d(X)
\cap \H^p(X)$ is the group of codimension $p$ cycles which are linear
combinations of cup-products of divisors. If $\d(X) = \H(X)$, then the
(1,1)-theorem implies the Hodge conjecture.

On the other hand, for an algebraic variety defined
over an algebraic number field, say $K \subset \Qbar$, one can consider
$\ell$-adic \'etale
cohomology $H\,\dot{ }_{\!\!\et}(X_\Qbar, \Ql).$
The Galois group $\Gal(\Qbar/K)$ acts continuously on
$H\,\dot{}_{\!\!\et}(X_\Qbar, \Ql)$. 
If $F$ is a finite extension of $K$, then the open subgroup $\Gal(\Qbar/F)$ of 
$\Gal(\Qbar/K)$ acts by the $p^{th}$ power of the inverse of the cyclotomic 
character $\chi_\ell$ on the cohomology classes of $F$-rational codimension $p$
algebraic cycles. If for an arbitrary $\ell$-adic Galois representation 
$W$ and an integer $n \in \Z$ $W(n) := W \otimes \chi_\ell^n$ denotes the 
$n^{th}$ Tate twist of $W$, then the these cohomology classes are in 
$H^{2p}_{\!\et}(X_\Qbar, \Ql)(p)^{\Gal(\Qbar/F)}$. By a {\it codimension 
$p$ Tate cycle} we mean a cohomology class in $H^{2p}_{\!\et}(X_\Qbar, \Ql)(p)$
fixed by an open subgroup of $\Gal(\Qbar/K)$. Since 
$H\,\dot{ }_{\!\!\et}(X_\Qbar, \Ql)$ is a finite dimensional $\ell$-adic 
representation, we can find the largest open subgroup of the Galois group fixing
all the Tate cycles (hence these cycles are all defined over a large enough
number field).
Let $\g$ be the Lie algebra of the image
of the Galois group in $\End_{\Q_l}(H\,\dot{ }_{\!\!\et}(X_\Qbar, \Ql))$. 
Then codimension $p$ Tate cycles $\T^p(X)$ are, by definition, the
$\g$-invariants $H^{2p}_{\!\et}(X_\Qbar, \Ql)(p)^\g$.
Tate [T 0] conjectured that {\it all the Tate cycles are algebraic}.
We denote $\T:=\underset {p} \to \oplus \, \T^p$ the ring (under the
cup-product) of the Tate cycles.

{}From now on, we restrict ourselves to the case of {\it
abelian varieties} defined over number fields. On the one
hand, this case is more concrete, and some progress has been made;
on the other, it has important arithmetical applications.

Although not known in general, the analog of the (1,1)-theorem for the Tate
cycles of codimension 1 for \avs\ has been proved by Faltings [F]. Hence, as
above, we can conclude that the Tate conjecture holds for an \av\ $A$
satisfying $\T(A) = \d_{\ell}(A)$, where $\d_{\ell}(A)$ is the ring of
\'etale cohomology classes generated by divisors.

It is known that {\it generically}, but not always, the Hodge (resp. the Tate)
cycles are all generated by divisors [Ma], [Ab~1].

The first (counter)example due to Mumford (cf.\ [Po]) features a CM abelian
4-fold.
Weil [W] has shown that the essential feature of Mumford's example causing
$\H=\d$ to fail is an action in a special way of a quadratic imaginary field
$\k$  on an abelian variety. Namely, consider a family of \avs\ of even
dimension, say $2d$, whose endomorphism algebra contains such a field $\k$ with
the signature of the $\k$-action $(d, d)$. Generically for such a family, the
ring
of Hodge cycles is generated by divisors together with an exceptional
(non-divisorial) cycle of codimension $d$.

Recently, C. Schoen proved the Hodge conjecture for one family of abelian
4-folds of
Weil type (with an action of $\Q(\mu_3)).\ $

Our {\it initial} motivation was to answer a question of Tate on 
whether the Tate conjecture holds for this family (cf.\ [T~2], p.~82). 
The affirmative answer was obtained independently by Moonen-Zarhin 
(cf.\ \cite{MZ}).

In general, both conjectures seem to be very difficult in codimensions $> 1$.

The existence of the comparison isomorphisms
between the $\ell$-adic and singular cohomology theories carrying
algebraic cycles in one theory to another suggests that the Hodge and the Tate
conjectures describe essentially the same object. So, it is
natural to ask if the two conjectures are equivalent in some sense.
The precise statement in the case of abelian varieties constitutes
the {\ffi Mumford--Tate conjecture}, which we denote by {\bf MT}.
It asserts that {\ffi the Hodge and the Tate conjectures are equivalent for an
abelian variety and all its self-products}.

Concretely, for an \av\ $A$ over a \nf\ $K$, there exists a connected reductive
algebraic subgroup
$\Hg(A)$ of $\GL(V)$ defined over $\Q,\ V =$
\linebreak
$H^1(A_\C, \Q)$
(resp.\ a connected reductive algebraic subgroup $G_\ell(A)$ of $\GL(V_\ell)$
defined over $\Ql,\ V_\ell = H_\et^1(A_\Qbar, \Ql)$ for some prime number $\ell
\! \in \! \Z$),
\s.t. the Hodge (resp.\ the Tate) cycles of codimension $p$ are obtained as
invariants in $H^{2p}(A, \Q) \iso \wedge^{2p} V$
(resp.\ in $H_{\et}^{2p}(A_{\Qbar}, \Ql) \iso \wedge^{2p} V_\ell$)
of $\h = \Lie(\Hg(A))$ (resp.\ $\g_\ell = \Lie(G_\ell(A))$). Because of the
comparison isomorphism $V_\ell \iso V \otimes_\Q \Ql$ between the two
cohomology theories, $\Hg_\ell(A) := \Hg(A) \otimes_\Q \Ql$ acts on
$V_\ell$.
Let $\h_\ell := \Lie(\Hg_\ell) = \h \otimes \Ql$.
MT asserts that $\g_\ell = \h_\ell$.
(Note that $\g_\ell(A)$ is {\it not} the Lie algebra of the image of
$\Gal(\Qbar/K)$ in $\GL(V_\ell)$, but the intersection of the Lie algebra of
this image with $\sl(V_\ell)$.
That is why we do not Tate-twist the \'etale cohomology group.)

Deligne, Piatetskii-Shapiro and Borovoi proved a ``half" of MT, viz.,
$\h_\ell(A) \supseteq \g_\ell(A)$. Hence the Tate conjecture implies the Hodge
conjecture.

MT for \avs\ of CM-type is a consequence of the results of Shimura and Taniyama
(cf.\ [ShT], [Po]).
This must have been the motivating factor behind Mumford--Tate.

MT has been proved in a few (non-CM) cases by imposing restrictions on the
size of the endomorphism algebra and adding some divisibility conditions on the
dimension of \avs, cf.\ [S 0], [C 0]. In [Z~4], [Z~5], [LZ] 
MT was verified for \avs\ satisfying certain conditions
on the Galois action at a prime of good reduction. In all these cases 
the Tate cycles are generated by divisors, hence the Tate conjecture holds. 
The result now follows from the ``known half" of MT.

The main thrust of this work is to show that under suitable bad reduction
conditions we
can control the image of Galois; in particular, MT holds for a class of abelian
varieties, including some Weil-type abelian varieties (for which the Tate 
conjecture is {\it not} known).

Note that if $A$\ is an absolutely simple \av, $e =  (\End^{\circ}(A) :
\Q)$ the degree over $\Q$ of its endomorphism algebra, and $A$ has bad \red\ at
some prime $\wp$,\ then $e$ divides the dimension of the toric part of the
reduction.

The following is the main result of this work (see Theorem 6.4)
\proclaim
{Main Theorem} Let $A$\ be an absolutely simple \av, $\End^{\circ}(A) =
\Bbbk$ : imaginary
quadratic field, $g = \dim(A)$. Assume $A$ has bad semi-stable reduction 
at some prime $\wp$, with the dimension of the
toric part of the reduction equal to $2r$, and $\gcd(r, g) = 1$, and $(r,
g) \neq (15, 56)$ or $(m-1, \frac {m(m+1)}2)$. 
Then MT holds.
\endproclaim

Roughly speaking, the idea is the following. If $A$ is an abelian variety 
with bad semi-stable reduction
at some prime $\wp$ (of its field of definition), then the action of the
inertia at $\p$ on on the $\ell$-adic ($\p \nmid \ell$) Tate module of 
$A$ is unipotent of ``rank" 
equal to the dimension of the toric part of the identity component of 
the the special fiber of the N\`eron model of $A$ at $\p$ (=: toric rank). 
If $A$ satisfies the 
conditions on the ``size" of the endomorphism algebra and the toric rank
imposed above, then the rank of the unipotents (in the inertia image) is
prime to the dimension of the Galois representation, which is a 
sufficiently restrictive condition, given our knowledge of the possible
Galois representations arising in this situation.

Note that people have looked at the special elements in the monodromy 
action before (cf.\ 6.8).

As mentioned above, MT is known to hold for CM \avs. It is also known (e.g.,
[ST]) that such \avs\ have good \red\ at all primes, after possibly a finite
base change. 
But the set of abelian varieties with (potentially) good reduction 
{\it everywhere} is ``small" in a corresponding moduli space, i.e., it is a 
very rare occasion that an abelian variety has everywhere (potentially)
good reduction.
Indeed, such \avs\ correspond to ``integral" points of the moduli space 
(cf.\ \loc.cit., Remark (1), p.~498)
and as ``sparse" as integers in a number field.
So, ``most" of the \avs\ do have bad reduction ``somewhere."
We have reasons to believe that 
\avs\ with {\it minimal} bad reduction (e.g., the case $r = 1$ of Theorem A)
are the ``most typical" (cf.\ [L]).

\proclaim\nofrills{\rm \quad 
Along the way we established various other results. They are: }
\roster
\item"$\bullet$" 
If $A$ is an abelian variety with $\End^{\circ}(A) = \Q$ 
and the dimension of the toric part of its reduction 
is either $2$ or prime to $2\dim(A)$, then MT holds (Theorems 6.5, 6.6).

\item"$\bullet$" MT holds for some abelian 4-folds $A$ with
$\End^{\circ}(A) = \Q$ (Theorem 5.2).

\item"$\bullet$" For some abelian varieties, either MT or the Hodge conjecture
holds (Theorem 7.1, Remark 7.2(4)).
\endroster

\endproclaim

The question of existence of \avs\ considered here and the ``size" of the set
of such varieties in the corresponding moduli spaces is addressed in [L].

\medskip
This work started with 4-dimensional case, part of which was 
independently obtained by Moonen-Zarhin. Since this has been published 
(\cite{MZ}), we do not treat this case here (see, however, section 2).
\subhead{\it Acknowledgment}\endsubhead
I greatly benefited from conversations with many people over the 
period of time I was working on this paper.
I would like to thank my advisor, Dinakar Ramakrishnan, for encouragement 
and guidance throughout the research and preparation of my 
PhD thesis (at Caltech) in which most of these results appeared.
I am especially grateful to Don Blasius, Haruzo Hida and
Michael Aschbacher for numerous helpful remarks and enlightening comments. 

This work could not have been possible without generous
help of Fedor Bogomolov. I can hardly thank him enough.

\head
{\ff I.  Preliminaries, first applications}
\endhead
\head
{{\ff 0. Basics and Notation}}
\endhead

\subhead
{0.0}
\endsubhead
Let $A\ $be a {\it simple} (:=absolutely simple) abelian variety defined over
some number field,
say, $K \embed \Qbar$ with a fixed embedding, $D := \End^{\circ}(A) :=
\End_\Qbar^{\circ}(A) \otimes_\Z \Q,\ V:=H_1(A(\C),\Q)$. Then
$D \hookrightarrow \End_\Q(V).$

We recall the Albert classification of the possible {\it types} of the
endomorphism algebras $D$
(cf.\ [MAV, $\S$20]):
\roster
\item"{\bf I}\ :" $D$ is a totally real field;
\item"{\bf II}\ :" $D$ is a {\it indefinite} quaternion algebra over a totally real
field $F$,
i.e., $D \otimes_F \R$ is a sum of $f := (F:\Q)$ copies of M$_2(\R)$;
\item"{\bf III}\ :" $D$ is an {\it definite} quaternion algebra over a totally real
field $F$,
i.e., $D \otimes_F \R$ is a sum of $f := (F:\Q)$ copies of the
Hamiltonian quaternions ${\Bbb H}$;
\item"{\bf IV}\ :" $D$ is a division algebra over a $CM$-field;
\endroster

By abuse of language we say that $A$ is ``of type I (II, ...)" if its $D$
is of this type.

\subhead
{0.1}
\endsubhead
Recall that $V_\R := V \otimes_\Q \R\ $is given a complex structure induced by
the
natural isomorphism between $V_\R\ $and the universal covering space of $A(\C)\
$(cf.\ [MAV]). Therefore we obtain a homomorphism of algebraic groups,
$$\varphi : T \to \GL(V),$$
defined over $\R,\ $where $T\ $is the compact one-dimensional torus over $\R,\
$i.e., $T_\R =
\{z \in \C\ |\ |z| = 1\},\ $by the formula
$$
\aligned
\varphi(e^{i \theta}) = &\text{ the element of} \GL(V)\text{, which is
multiplication by }e^{i \theta}\\
&\text{ in the complex structure on }V_\R.
\endaligned
$$
Note that there is a non-degenerate skew symmetric (Riemann) form $\Theta : V
\times V \to \Q$
and that $\varphi$ satisfies the Riemann conditions (cf.\ [M~1]):
\roster
\item"1." $\varphi (T) \subseteq \text{Sp}(V, \Theta),$
\item"2." $\Theta(v, \varphi(i) \cdot v) > 0,\ \forall v \in V,\ x \neq 0.$
\endroster

\definition
{Definition}
The {\it Hodge\ group}\ $\Hg(A)$ of $A$ is the
smallest algebraic subgroup of Sp$(V) :=\ $Sp$(V, \Theta)$ defined over $\Q$
which after extension of scalars to $\R$ contains the image of $\varphi$.
\enddefinition

\subhead
{0.1.1}
\endsubhead
For the purpose of completeness and further reference, we give the following
reformulation of the above construction and definition. The reference for what
follows is [D~3, $\S\ $3].

Since $A(\C)\ $is a compact smooth K\"ahler manifold, $V_\C := H_1 (A(\C), \C)$
admits a
Hodge decomposition
$$H_1 (A(\C), \C) = H_{-1,0} (A) \oplus H_{0,-1} (A).$$
Thus we obtain a homomorphism
$$\mu: {\Gm}_{,\C} \to \GL(V)_\C$$
by defining $\mu(z),\ \forall z \in \C^{\sssize \times},\ $to be the
automorphism of $V_\C\ $which
is multiplication by $z\ $on $H_{-1,0} (A)$ and by the identity on $H_{0,-1}(A)$.

\defn The {\it Mumford--Tate group} $\M(A)$ of $A$ is the smallest algebraic
subgroup of $\GL(V)$ defined over $\Q$ which after extension of scalars to
$\C\ $contains the
image of $\mu.$
\enddefn

Clearly, over $\C,\ \M(A)$ is the subgroup of $\GL(V)_\C$ generated by the
conjugates
$^\sigma \! \mu,\ \forall \sigma \in\ \Aut(\C).$

\defn The {\it Hodge\ group}\ $\Hg(A)$ of $A$ (or the {\it special}\
Mumford--Tate
group of $A$) is the connected component of the identity of the intersection
$\M(A)$ $\cap$ SL$(V)$ in $\GL(V)$.
\enddefn

\rems 1.  The construction of $\M(A)$ furnishes it with a canonical character
$\nu : \M(A) \to \Gm$ defined over $\Q$ and characterized by the condition 
$\nu \circ \mu =\ $id$_\Gm$. Then $\Hg(A) = \Ker(\nu)$. This is the 
reason why we use the Hodge group instead of the Mumford--Tate group 
(cf.\ also 0.3).

2.  One can easily show that the two definitions of $\Hg(A)$ are
equivalent.
\endrem

\subhead
{0.1.2}
\endsubhead
The following theorem lists the properties of $\Hg(A)$ 
(cf. \cite{M 1, $\S$ 2, Theorem (i)}, \cite{Ta, Lemma 1.4}).

\thm
\roster
\item"1."  $\Hg(A)$ is a connected
reductive group.
\item"2."  $D\ (= \End^\circ (A)) =
\End_{\Hg(A)}(V)\ = \End_{\h}(A)$,
where $\h:={\Lie}(\Hg(A))$.
\item"3."  $\Hg(A)$ is \ss\ for an abelian variety $A$ of type I, II or III.
\item"4."  $\Hg(A^a \times B^b) \iso \Hg(A \times B)$ for 
any \avs\ $A, B$ and $a,b \in \N$.
\endroster
\endst

\rems 1.  Part 2 of the theorem implies that $A$ is simple if and only if
$V$ is $\h$-simple, if and only if (Schur's lemma) $D$ is a division algebra.

2.  One can refine part 3, cf.\ 0.7.1.
\endrem

\subhead
{0.1.3}
\endsubhead
Recall that the {\it Hodge\ classes}\ of $A\ $are classes of type $(p, p)\ $in
the
Hodge decomposition of homology of $A.$

The {\it Hodge\ conjecture}\ states that all the Hodge classes are algebraic.

\subhead
{0.1.4}
\endsubhead
By the K\"unneth formula $H_{\text{\bf .}}(A) = \overset {\text{\bf .}} \to \wedge \, H_1(A)\ $(cf.\ [MAV]), hence $\Hg(A)$ acts on 
$H_{\text{\bf .}}(A)$. One can show (cf.\ [M~1]) that the Hodge 
classes of $A\ $are exactly those classes in $H_{\text{\bf .}}(A)$
that are fixed
by $\Hg(A)$. In fact, the Hodge group has the following {\it characteristic}
property (\loc.cit., $\S$ 2, Corollary).

\thm
The Hodge group $\Hg(A)$ is the largest (reductive)
subgroup of $\GL(V)$ fixing all the Hodge classes of $A^s,\ s \geq 1.$
\endst

By the K\"unneth formula $H_2(A^s) = \overset s \to {\underset {i=1} \to
\oplus}
H_2(A).\ $Hence, in the view of the previous theorem, the Lefschetz
(1,1)-theorem for abelian
varieties takes the following form.

\thm Let $s\in\N,\ sV:= V \oplus ...
\oplus V (s$ times$)$,
then the $\h$-invariants  $(\underset{\Q}\to{\overset 2\to \wedge}\, sV)^\h\ $
is exactly the $(\Q$-span of homological classes of$)$ divisors on $A^s = A
\times ... \times A\ (s$ times$)$.
\endst
\subhead
{0.1.5}
\endsubhead
Following Ribet and Murty (cf.\ \cite{Ri}, \cite{Mu}) we make the following

\defn
The {\it Lefschetz\ group}\ $L(A)$ of $A$ is the
connected component of the identity of the centralizer of $\End^{\circ}(A)$ in
Sp($V, \Theta)$ (inside $\End_\Q(V),\ \Theta$ is a polarization, cf.\ 0.1).
\enddefn

The following are the main results about the Lefschetz group (cf.\ \cite{Mu},
also \cite{Ri}, \cite{H 1}).

\thm
\roster
\item"0.(i)"  $L(A)\ $is a connected reductive algebraic group
defined over $\Q.$
\item"(ii)"  $\Hg(A) \subseteq L(A).\ $
\item"(iii)"  $L(A)\ $is semi-simple for $A\ $of type I, II or III; moreover, it is
symplectic for $A$ of type I and II, orthogonal for type III.
\item"(iv)"  $L(A^{n_1}_1 \times ... \times A^{n_s}_s) = L(A_1) \times ... \times
L(A_s).$

\item"1.   "  All the Hodge classes on $A^s\ $are divisorial if and 
only if $\Hg(A) = L(A)$ and $A$ is not of type III.
\item"2.   "  If $A$ is of type III, then 
it has a non-divisorial Hodge class.
\endroster
\endst

Let $\l\ $:= $\Lie(L(A)$), then
\roster
\item"0.(ii)$'$" \ \   $\h \embed \l \embed \sp(V),\ \h^{ss} \embed \l^{ss},\
C_{\h} \embed C_{\l}$.
\endroster

Here $C_?$ is the center of $?$, and $?^{ss}$ is the \ss\ part of $?$.

\subhead
{0.2}
\endsubhead
Let $\k \hookrightarrow D$ be an imaginary quadratic field,
$\Gal(\k/\Q) = \{\sigma, \rho\}, \rho\ (= {\sigma}^{2})$ is the fixed
(identity) embedding $\k \embed \Qbar$. In this case $V_\R := V \otimes \R$ 
has two complex structures. One is given by the isomorphism 
$V_\R = \Lie(A(\C))$,\ (cf.\ 0.1), and the other by the action of 
$\k \otimes_\Q \R (\simeq \C)$. Hence the splitting $V_\R = V^\sigma \oplus
V^\rho$ ($\k$ acts by $\sigma(\k)$ on $V^\sigma$ and by $\rho(\k)$ on $V^\rho$).
The two complex structures coincide on one of the subspaces, say $V^\rho\ $, and
conjugate on the other, $V^\sigma$. If $m_\sigma = \dim_\C(V^\sigma),\
m_\rho = \dim_\C(V^\rho)$, then
$(m_\sigma, m_\rho)\ $is the {\it signature}\ of the $\k$-action; $m_\sigma +
m_\rho = g = \dim_\C(V_\R) = \dim(A)$.

\subhead
{0.2.1}
\endsubhead
Recall that the {\it Rosati\ involution}\ is the involution
on $D = \End^{\circ}(A)$
induced by the Riemann form (cf.\ 0.1). The Rosati involution is {\it
positive}, consequently, the field it fixes is totally real (cf.\ [MAV]).

In the case $\k \subseteq D$ we always assume that the {\it Rosati involution
preserves} $\k$. The
positivity of the involution implies that it acts on $\k$ non-trivially. Hence
this action coincides
with (the complex conjugation) $\sigma.\ $

\subhead
{0.2.2}
\endsubhead
Since $\h$ and $\l$ centralize $D$ (cf.\ 0.1.2  and 0.1.5)
$$
\h \embed \l \embed \sp(V)^\k \embed \sp(V),
$$
where $\sp(V)^\k\ $is the centralizer of $\k\ $in $\sp(V)$. By [D~3, Lemma 4.6]
$$
\sp(V)^\k = \u(V),
$$
the Lie algebra of the unitary group of a $\k$-Hermitian form on $V\ $viewed as
the $\k$-vector space (cf.\ 0.2). Extending scalars to $\k\ $we get
 $$\h_\k \embed \l_\k \embed \u(V) \times \u(V)^\sigma \embed \sp(V_\k),
\tag 0.2.2.1$$
where $\h_\k:=\h \otimes_\Q \k,\ \l_\k:=\l \otimes_\Q \k,\ V_\k := V \otimes_\Q
\k = V \oplus U,\ U\ $is
the same $V,\ $but with the conjugate $\k$-vector space structure.

The $\h$-invariant $\k$-Hermitian form referred to above
is a non-degenerate
element of $\check V \otimes \check U$ (cf.\ [D~3, Lemma 4.6]), hence the
isomorphism
$$U \iso \check V$$
of $\h$-modules. Clearly the projection of $\h_\k\ $to $\u(V)\ $is $\h,\ $thus
we can rewrite
(0.2.2.1) as
 $$\h \embed \l \embed \u(V) \overset{\, \Delta} \to {\embed} \sp(V \oplus
\check V),\      g
\overset{\Delta} \to {\mapsto}  \pmatrix g & 0 \\ 0 & ^tg^{-1}
\endpmatrix. \tag 0.2.2.2$$
\remark{{\rm (}Note} 
The embeddings above are considered over $\k$.)
\endnote

\subhead
{0.2.3}
\endsubhead
From this we get:
 $$\overline \h \embed \overline \l \embed  \gl(W) \overset{\Delta} \to
{\embed} \sp(W \oplus
\check W),$$
where $\overline \h := \h \otimes_\Q \overline \Q,\ \overline \l := \l
\otimes_\Q \overline \Q,\ W :=
V \otimes_{\k,\rho} \overline \Q,\
\check W\ $: dual of $W\ (= W \otimes_{\k,\sigma} \overline \Q).\ $

\subhead
{0.3}
\endsubhead
Let $V_\ell := T_\ell(A) \otimes_\Zl \Ql \iso V \otimes_\Q
\Ql,\ $where $T_\ell(A)\ $is
the $\ell$-adic Tate module ($=~H_1^\et(A\times_K \Qbar,\Zl) = \underset n \to
\varprojlim \,
{}_{\ell^n} \! A(\Qbar),\ $where ${}_{\ell^n} \! A(\Qbar) =\ $kernel of
multiplication by $\ell^n : A(\Qbar) \to A(\Qbar)),\ V = H_1(A(\C), \Q)\ $as
above. By abuse of language we call $V_\ell$ Tate module too. Let 
$G_\ell$ be the image of $\Gal(\Qbar/K)$ in
$\End_\Ql(V_\ell)$, where $K$ is the base
field of $A,\ \G_\ell := \Lie(G_\ell)$.
It is known (cf.\ [Bo]) that $\G_\ell$ is algebraic and 
$\Ql \cdot 1_{V_\ell} \subseteq \G_\ell \subset \g\sp(V_\ell)$.
Let $\g_\ell := \G_\ell \cap \sl(V_\ell) \subset \sp(V_\ell)$. Then
$C_{\G_\ell} =C_{\g_\ell} \oplus \Q_\ell \! \cdot \! 1_{V_\ell},\ \g_\ell^{ss}
= \G_\ell^{ss}.$

\rems 1. $\g_\ell\ $does not depend on finite extensions of $K\ $(cf.\ [S~2]).
\endrem

2. $V_\ell \iso V_\ell(r),\ \forall r \in \Z,\ $as $\g_\ell$-modules, but not
as $\G_\ell$-modules.

\subhead
{0.3.0}
\endsubhead
The {\it Tate\ conjecture}\ states that the {\it Tate\
cycles,}\ i.e., the Galois invariants $H_{\text{\bf .}}^\et(A_\Qbar,
\Ql)^{\g_\ell} = (\overset{\text{\bf .}}\to \wedge \, H_1^\et(A_\Qbar,
\Ql))^{\g_\ell}$, are algebraic.

\subhead
{0.3.1}
\endsubhead
Faltings (cf.\ \cite{F}) has proved the analogs of Mumford's
Theorem 0.1.2(1,2) and the
(1,1)-theorem (a special case of the Tate conjecture) for \avs.

\thm
\roster
\item"1. "  Let $s \in \N,\  (\underset{\,\,\Ql}\to{\overset 2\to \wedge}\,
sV_\ell)^{\g_\ell}\ $
is exactly the ($\Ql$-span of homological classes of) divisors of $A^s = A
\times ... \times A\  (s\
$times).
\item"2. "  $\End_{\g_\ell}(V_\ell) = D \otimes_\Q \Ql.$
\item"3. "  $\g_\ell$ is reductive.
\endroster
\endst

\subhead
{0.3.2}
\endsubhead
$\h_\ell := \h \otimes_\Q \Ql \embed \End_\Ql(V_\ell)$.
The known relation
between $\g_\ell\ $and $\h_\ell\ $is given by the following theorem.

\thM{Deligne \cite{D 3}, Piatetskii-Shapiro \cite{P-Sh}, Borovoi \cite{Bor}}

$\g_\ell \subseteq \h_\ell  \subset \End_\Ql(V_\ell).$
\endst

\subhead
{0.4}
\endsubhead
The {\it Mumford--Tate\ conjecture}\ (=: MT) states $\g_\ell =
\h_\ell.\ $Since $\g_\ell\
$and $\h_\ell\ $are reductive, it is the same as equivalence of the Hodge and
the Tate
conjectures for an abelian variety and all its self-products.

\subhead
{0.4.1}
\endsubhead
In order to prove MT it is enough to establish the
conjecture for one $\ell$ ([LP, Theorem 4.3]).

\subhead
{0.4.2}
\endsubhead
Moreover, it is enough to show $\overline \g_\ell =
\overline \h_\ell,\ $where
$\overline \g_\ell := \g_\ell \otimes_\Ql \Qlbar,\  \overline \h_\ell :=
\h_\ell \otimes_\Ql \Qlbar\
$([Z~2, $\S 5,\ $Key Lemma]).

\subhead
{0.4.3}
\endsubhead
The (1,1)-theorems imply $(\underset
{\,\,\Ql} \to {\overset 2 \to
\wedge} sV_\ell)^{\g_\ell} = (\underset {\,\,\Ql} \to {\overset 2 \to \wedge}
sV_\ell)^{\h_\ell}. $

\subhead
{0.4.4}
\endsubhead
Theorems 0.3.2, 0.1.2(2) and 0.3.1(2) imply $\g_\ell^{ss}
\subset \h_\ell^{ss},\ C_{\g_\ell}
\subset C_{\h_\ell}$.

{\ffi In fact}, $C_{\g_\ell} = C_{\h_\ell}.\ $This can be shown in a way similar
to the proof of MT for
CM abelian varieties (cf.\ [ShT], and also [D~1]).

So, {\ffi in order to prove MT one must show that 
$\overline\g_\ell^{ss} = \overline\h_\ell^{ss}$}.

To simplify notations we write sometimes $\gg$ for $\overline\g_\ell^{ss}$
and $\hh$ for $\overline\h_\ell^{ss}$.

\subhead {0.5} \endsubhead
Let $\a\ $be a \ss\ Lie algebra over an algebraically closed
field of characteristic 0,
and $\a = \a_1 \times ... \times \a_n\ $be the decomposition of $\a\ $into the
product of its simple
ideals. For any faithful \irrep\ $U$ of $\a,\ U$ decomposes as a tensor
product of \irreps\ $U_i$
of $\a_i$. Since $U$ is faithful, none of the $U_i$'s is trivial. Moreover,
if the \rep\ $U$ admits
a non-degenerate invariant bilinear form, then so does each $U_i.$

We say that the \rep\ is {\it minuscule}\ if the highest weight of each $U_i\
$is minuscule, see
[B, Ch.VIII, $\S7.3]$. The following is the list of minuscule weights, [B,
Ch.VIII,
$\S7.3\ $and Table 2]:

type $\bold{A_m}\ (m \geq 1):\ \varpi_1,\ \varpi_2,\ \dots\ ,\
\varpi_m;\  \dim(\varpi_s) =
\binom{m+1}s;$

type $\bold{B_m}\ (m \geq 2):\ \varpi_1;\  \dim(\varpi_1) = 2m+1;$

type $\bold{C_m}\ (m \geq 2):\ \varpi_1;\  \dim(\varpi_1) = 2m;$

type $\bold{D_m}\ (m \geq 3):\ \varpi_1, \varpi_{m-1}, \varpi_m;\
\dim(\varpi_1)=2m,$ \newline \hphantom{\indent type $D_m\ (m \geq 3):\ \ $}
$\dim(\varpi_{m-1}) = \dim(\varpi_m) = 2^{m-1};$

type $\bold{E_6}\ :\ \varpi_1,\ \varpi_6;\  \dim(\varpi_1) =
\dim(\varpi_6) = 27;$

type $\bold{E_7}\ :\ \varpi_7;\  \dim(\varpi_7)$ = 56;

there are {\it no} minuscule \reps\ for the types $\bold{E_8},\
\bold{F_4},\ \bold{G_2}$.

\subhead {0.5.1} \endsubhead
It is known that the representations of $\gg,\ \hh$ 
are minuscule (cf.\ [S~$\ast$], [D~2]). It is also known that $\gg$
is not exceptional, see [S~$\ast$, Theorem 7] (for the corresponding result for
$\hh$
see [D~2, Remarque 1.3.10(i)]).

\subhead {0.6.1} \endsubhead
Let again $\k \embed D\ $and $\k_\ell := \k \otimes_\Q \Ql.\
$Then, as in 0.2.2,
$$\g_\ell \subset \h_\ell \subset \sp(V_\ell)^{\k_\ell} \subset \sp(V_\ell).$$
If $l\ $splits in $\k,\ \lambda,\lambda'\ $being the primes of $\k\ $over
$l,\ \lambda' = \lambda^\sigma,\ $then $\k_\ell \iso \Ql \oplus \Ql,
V_\ell = V_\lambda \oplus V_{\lambda'}$, where $V_\lambda, V_{\lambda'}$ are
vector spaces over
$\k_\lambda \iso \Ql,\ \k_{\lambda'} \iso \Ql$ respectively, and
$$\g_\ell \subset \h_\ell \subset \gl(V_\lambda) \oplus \gl(V_{\lambda'})
\overset{\Delta} \to
{\subset} \sp(V_\lambda \oplus V_{\lambda'}).$$
Since $\lambda' = \lambda^\sigma$, as in 0.2.2 we conclude
$V_{\lambda'} \iso \check V_\lambda\ $and can rewrite the above sequence as
$$\g_\ell \subset \h_\ell \subset \gl(V_\lambda) \overset{\Delta} \to {\embed}
\sp(V_\lambda
\oplus \check V_\lambda).$$

\rem If $\ell\ $does {\it not}\ split in $\k,\ $then $\k_\ell\ $is a
field, $(\k_\ell : \Ql) = 2$, and
the rest is identical to 0.2.2.
\endrem

\subhead {0.6.2} \endsubhead
As in 0.2.2, by extending scalars to $\Qlbar\ $we get
$$\overline \g_\ell \subset \overline \h_\ell \subset \gl(W_\lambda)
\overset{\Delta} \to {\embed}
\sp(W_\lambda \oplus \check W_\lambda),$$
where $W_\lambda := V_\ell \otimes_{\k_\ell,\rho_\lambda} \Qlbar,\ \rho_\lambda
: \k_\ell \to
\k_\lambda\ $is the projection.

\rem For $\ell\ $non-split in $\k\ $, the same holds (cf.\ 0.2.2, 0.2.3).
\endrem

\subhead {0.7} \endsubhead
We will need the following simple facts. We assume that $D = \k$.

\Prop{0.7.1} The \reps\ of $\overline \g_\ell$ and $\overline \h_\ell$ 
are non-self-dual. \QED
\endst

\rems 1.  This is true for any \irr\ sub\rep\ of $W_\lambda$ for any
type IV \av\ (e.g., [Mu], [H]).

2.  If the \av\ is of type I (respectively II, respectively III), then the
\irr\ components are
symplectic (respectively symplectic, respectively orthogonal) (\loc.cit.).
\endrem

\Prop{0.7.2} $\overline \g_\ell$ and $\overline \h_\ell$ are \ss\
if and only if the signature of the $\k$-action is $(m,m).\ $
Further, if this is {\it not}\ the case, the centers $C_{\overline \g_\ell},
C_{\overline \h_\ell}\ $are 1-dimensional.
\endst

\pf This is essentially proved in [D~3], [W]. Let us,
however, briefly explain why this holds and fix notations.

Let $\mu : \Gm_{,\C} \to \GL(V_\R)$ be the cocharacter defining the Hodge
structure on
$V$, then the map $h : \Bbb S = R_{\C/\R}\Gm \to \GL(V_\C)$ is given by $h(z)
= \mu(z)$ on
$V_\R$ (cf.\ 0.1.1), $h(z) = {\overline \mu(z)}$ on $U_\R$, where $U_\R$ is
the ``same" $V_\R$ but with the
conjugate $\k \otimes \R$-action, $V_\C  = V_\R \oplus U_\R$ (cf.\ 0.2). If
the $\k$-signature is
$(m_\sigma,m_\rho)$ then $V_\R = V^\sigma_\R \oplus V^\rho_\R,\
\dim_\C(V^\sigma_\R)
= m_\sigma,\ \dim_\C(V^\rho_\R) = m_\rho\ (\k\ $acts by $\sigma(\k)$ on
$V^\sigma_\R$
and by $\rho(\k)\ $on $V^\rho_\R).\ $But $V_\R = H_{-1,0}$, hence the power
of $z$ by which
$\mu(z)$ acts on $V^\sigma_\R,\ V^\rho_\R$ is 1. We will call this power the
$\mu$-{\it weight}. Similarly, $V_\C = H_1(A(\C),\C) = V^\sigma_\C 
\oplus V^\rho_\C.\ $But also $V_\C =
H_{-1,0} \oplus H_{0,-1} (= V_\R \oplus U_\R)$ and these two decompositions
commute, since
the former is determined by $\k \subseteq D\ $and the Hodge group centralizes
$D$ in $\End_\Q(V)$. Hence we can write
$$\align
V_\C  &= V_\R \oplus U_\R \\
   &= (V^\sigma_\R \oplus V^\rho_\R) \oplus (U^\sigma_\R \oplus U^\rho_\R)\\
   &= (V^\sigma_\R \oplus U^\sigma_\R) \oplus (V^\rho_\R \oplus U^\rho_\R)\\
   &= V^\sigma_\C \oplus V^\rho_\C,
\endalign$$
where $U^\sigma_\R\ $(respectively $U^\rho_\R$) is the conjugate of $V^\rho_\R$
(respectively $V^\sigma_\R$). Thus 
\linebreak 
$\dim_\C(U^\sigma_\R)$ = $ m_\rho,\ \dim_\C(U^\rho_\R)
= m_\sigma,\ $so $\dim_\C(V^\sigma_\C) = m_\rho + m_\sigma = g =
\dim_\C(V^\rho_\C).$

The $\mu$-weight of $U_\R\ $is 0, hence the decomposition
$$V^\sigma_\C = V^\sigma_\R \oplus U^\sigma_\R $$
is according to $\mu$-weights 1, 0. (This exactly corresponds to the Hodge-Tate
 decomposition
of the $\lambda$-adic \rep below, cf.\ (1.0.$\ast\sim$1.1.$\ast$).) Now the
$\mu$-weight (respectively the Hodge type) of
$$\underset \C \to {\overset g \to \wedge} V^\sigma_\C =
\underset \C \to {\overset {m_\sigma} \to \wedge} V^\sigma_\R \otimes
\underset \C \to {\overset {m_\rho}  \to \wedge} U^\sigma_\R \subseteq
\overset g \to \wedge V^\sigma_\C$$
is $m_\sigma\ $(respectively $(-m_\sigma, -m_\rho)$). Hence $\underset \k \to
{\overset g \to
\wedge} V\ $is a Hodge cycle (i.e., of Hodge type
$(-\frac g2, -\frac g2)$) if and only if $m_\sigma = m_\rho.\ $Hence it is fixed
by $\h\ $if and only if $m_\sigma =
m_\rho,\ $i.e., $\h \subset \sl(V) \cap \u(V)\ $only in this case (here $V\ $is
considered as a
$\k$-vector space). In other words, the center $C_\h\ $of $\h\ $kills the
determinant $\underset \k \to
{\overset g \to \wedge} V\ $(and hence $\neq\ \{0\}$) if and only if $m_\sigma
\neq m_\rho.\ $Now, since
$\overline V := V \otimes \Qbar = W \oplus \check W, \ W$ is \irr,
$\overline \h \subset \gl(W) \overset {\Delta} \to {\to} \sp(\overline V)$.
Summarizing,
$$\align
\overline \h = \overline  \h^{ss}
&\subset \sl(W) \overset {\Delta} \to {\to} \sp(W \oplus \check W)\ \text{  if}\
m_\sigma = m_\rho,\\
\overline \h = \overline \h^{ss} \oplus \Qbar
&\subset \gl(W) \overset {\Delta} \to {\to} \sp(W \oplus \check W)\ \text{  if}\
m_\sigma \neq
m_\rho.
\endalign$$
Using 0.4.4 we conclude
$$\align
&\overline \g_\ell \subset \overline \h_\ell \subset \sl(W_\lambda) \overset
{\Delta} \to {\to}
\sp(W_\lambda \oplus \check W_\lambda),\\
&\overline \g_\ell = \overline \g^{ss}_\ell, \overline \h_\ell = \overline
\h^{ss}_\ell,\  \text{if}\ m_\sigma = m_\rho,\\
\\
&\overline \g_\ell \subset \overline \h_\ell \subset \gl(W_\lambda) \overset
{\Delta} \to {\to}
\sp(W_\lambda \oplus \check W_\lambda),\\
&C_{\overline \g_\ell} = C_{\overline \h_\ell} = \Qlbar,\  \text{if}\ m_\sigma
\neq m_\rho. \qed
\endalign$$
\endpf

\subhead {0.7.3} \endsubhead
As one can see from the proof, if $\k \subseteq
\End^{\circ}(A)\ $(not necessarily
equal), and $m_\sigma \neq m_\rho,\ $the center $C_\h\ $ of $\h\ $must kill the
determinant
$\underset \k \to {\overset g \to \wedge} V,\ $hence the center is non-trivial.
However, if $\k \neq
\End^{\circ}(A),\ $then the center can be non-trivial even if $m_\sigma =
m_\rho\ $(e.g., CM case).

But even if $\k \varsubsetneq \End^\circ(A)$, as follows from the proof,
$\overline \g_\ell \subset \overline \h_\ell \subset \sl(W_\lambda)$ if and
only if $m_\sigma = m_\rho$. In this form the result (using [D 3, Proposition
4.4]) can be generalized to the case of an arbitrary CM-field $E \embed
D$ (cf.\ [MZ, Lemma 2.8]).

Note that since $\g_\ell$ and $\h_\ell\ $are semi-simple for abelian 
varieties of types II or III (cf.\ 0.1.2(3)), if
$\k \embed \End^{\circ}(A)$ is Rosati-stable, $A$ of type II or III, then
the signature of the $\k$-action is
necessarily $(m,m)$.

\definition
{Definition}
If  the signature of the $\k$-action on an \av\ is $(m,m)$, we call
such an \av\ a {\it Weil\ type} \av\ (cf.\ [W]).
\enddefinition

\head
{\ff 1. On the Hodge-Tate decomposition}
\endhead

\subhead {1.0} \endsubhead
We recall here certain basic facts on the Hodge-Tate
decomposition and then give
some applications. The classical/standard reference is [T 1].

\subhead {1.0.1} \endsubhead
According to Tate and Raynaud, $\overline V_\ell := V_\ell
\otimes \Cl =
H_1((A \otimes_K \overline K)_{\et},\Ql) \otimes_\Ql \Cl,\ \Cl$ is a completion
of $\Qlbar$, admits a decomposition
$$
\overline V_\ell = \overline V_\ell (0) \oplus \overline V_\ell (1),
$$
where $\overline V_\ell (i) := \overline V_\ell^{(i)} \otimes_\Ql \Cl,\ i=1,2.\
$The $\Ql$-subspaces
(but {\it not}\ $\Cl$-subspaces) $\overline V_\ell^{(i)}$'s of  $\overline
V_\ell\ $are defined as
follows:
$$\overline V_\ell^{(i)} := \{ v \in \overline V_\ell \ |\ v^\sigma =
\chi_\ell(\sigma)^i \! \cdot \! v,\
\forall \sigma \in {\Cal I}
\},\ i=1,2,$$
where
${\Cal I}$ is the absolute inertia group at a prime of $K$ (=base field of $A$)
over $\ell$
(cf.\ 0.3 for why we take the inertia instead of the whole decomposition group),
$\chi_\ell$ is the cyclotomic character.
Recall that the Galois action is continous and semi-linear on
$\overline V_\ell\ $(see [S~2, 1.2]), and, clearly, the $\overline
V_\ell^{(i)}$'s are Galois
submodules of $\overline V_\ell.\ $The Galois
action on $\overline V_\ell (i)$ is by the formula
$$
(v \otimes c)^\sigma := v^\sigma \otimes c^\sigma,\ \forall v \in \overline
V_\ell^{(i)},\ \forall c \in \Cl,\ \forall \sigma \in {\Cal I},
$$
extended by linearity.

\subhead {1.0.2} \endsubhead
According to S. Sen ([Se, Section 4, Theorem 1]), to the
\HTd\ on $\overline V_\ell\
$one can associate a cocharacter
$$
\phi: \Gm_{,\Cl} \to \GL(V_\ell)_\Cl,
$$
by defining $\phi(z),\ \forall z \in \C^{\sssize \times}_\ell,\ $to be the
automorphism of $\overline V_\ell\
$which is multiplication by $z\ $on $\overline V_\ell(1)\ $and by the identity
on $\overline
V_\ell(0).\ $This association is made in such a manner that the algebraic
envelope $\tilde G_\ell$ of the Galois image (cf.\ 0.3) turns out to be the
smallest algebraic group defined over $\Ql\ $which after \ext\ of scalars to
$\Cl\ $contains the
image~of~$\phi.$

\rem This cocharacter $\phi$ is completely analogous to the
cocharacter $\mu$ associated to the \Hd\ on $V_\C,\ \tilde G_\ell$ is the
analog of the Mumford--Tate group $M(A)$ and
$\g_\ell$ is the analog of $\h$, see 0.1.1.
\endrem

\subhead {1.0.3} \endsubhead
Before proceeding, recall that for a
$\Gal(\Qbar_\ell/\Ql)$-module $X,\ $the {\it Tate\
twist}\ $X(1)\ $of $X\ $is defined to be $X \otimes_\Ql \Ql(1)\ $with the
Galois structure of a
tensor product of Galois modules (as in 1.0.1). Here $\Ql(1)\ $is the {\it
Tate\ module:}
$$
\Ql(1) := (\underset n \to \varprojlim \, \mu_{\ell^n}) \otimes_\Zl \Ql,\
\mu_{\ell^n} = \{ \zeta \in
\Qlbar\ |\ \zeta^{\ell^n} = 1\},
$$
with the natural $\Gal(\Qbar_\ell/\Ql)$-action by $\chi_\ell$:
$$
\zeta^\sigma = \zeta^{\chi_\ell(\sigma)},\ \forall \sigma \in
\Gal(\Qbar_\ell/\Ql),\ \zeta \in
\mu_{\ell^n}, \text{ for some}\ n.
$$

\subhead {1.0.4} \endsubhead
The Hodge-Tate decomposition of $\overline V_\ell\ $can be
rewritten in the
following explicit form ([T~1, $\S\ 4,\ $Corollary 2], see also the Remark
following that Corollary):
$$
\overline V_\ell = \Lie (A\spcheck_\Cl)\spcheck \oplus \Lie (A_\Cl)(1),
$$
where $\Lie(A\spcheck_\Cl)\spcheck\ $is the cotangent space of the dual \av\
$A\spcheck_\Cl\
$at its origin and $\Lie(A_\Cl)(1)\ $is the tangent space of $A_\Cl\ $at its
origin Tate-twisted by
$\chi_\ell.$

\subhead {1.0.5} \endsubhead
On the other hand, we have the Hodge decomposition on
$V_\C=H_1(A(\C),\C)$:
$$
\align
H_1(A(\C), \C) &= H_1(A_\C,\O_{A_\C}) \oplus H_0(A_\C, \Omega^1_{A_\C}) \\
&=  \Lie(A\spcheck_\C)\spcheck \oplus \Lie(A_\C),
\intertext{or, in our notation,}
V_\C &= U_\R \oplus V_\R,
\endalign
$$
see the proof of 0.7.2.

\subhead {1.0.6} \endsubhead
Fix an isomorphism $\Cl \iso \C.\ $Then the comparison
isomorphism
$$c: H_1^\et(A_{\overline K},\Ql) \otimes_\Ql \Cl \iso H_1(A(\C),\C)$$
or, in our notation,
$$c: V_\ell \otimes_\Ql \Cl \iso V \otimes_\Q \C,$$
provides isomorphisms
$$c: V_\R \to \overline V_\ell(1),\ c: U_\R \to \overline V_\ell(0). $$
Note, that $\dim_\C(V_\R) = g = \dim_\Cl(\overline V_\ell(1)),\
\dim_\C(U_\R) = g
= \dim_\Cl(\overline V_\ell(0)).$

\subhead {1.0.7} \endsubhead
On the other hand, both homology groups admit decompositions
according to the
action of $\k \subseteq \End^{\circ}(A):$
$$\overline V_\ell := \overline V_\lambda \oplus \overline V_{\lambda'}$$
where $\overline V_\lambda := W_\lambda \otimes_\Qlbar \Cl
= V_\ell \otimes_{\k_\ell,\rho_{\lambda}} \Cl
,\ \overline V_{\lambda'} := W_{\lambda'} \otimes_\Qlbar \Cl
= V_\ell \otimes_{\k_\ell,\rho_{\lambda'}} \Cl $, (cf.\ 0.6.2) and
$$V_\C := V_\C^\sigma \oplus V_\C^\rho.$$
These two types of splittings commute, and consequently $V_\C^\sigma,
V_\C^\rho\ $admit the
Hodge decomposition, and $\overline V_\lambda,\ \overline V_{\lambda'}\ $admit
the Hodge-Tate
decomposition. The map $c\ $ respects these splittings, hence maps either
$V_\C^\sigma\ $to $
\overline V_\lambda,\ $or $V_\C^\rho\ $to $ \overline V_\lambda.$

\subhead {1.0.8} \endsubhead
Let's {\it assume}\ that
$$c: V_\C^\sigma \overset \sim \to \to \overline V_\lambda,$$
hence $\dim_\C(V_\C^\sigma) = g = \dim_\Cl(\overline V_\lambda).\ $

As a result we conclude
$$c: V_\R^\sigma \overset \sim \to \to \overline V_\lambda(1),$$
hence $\dim_\Cl(\overline V_\lambda(1)) = m_\sigma,\
\dim_\Cl(\overline
V_\lambda(0)) = m_\rho.$

\rem  If $c: V_\C^\rho \to \overline V_\lambda,\ $then $m_\sigma\ $and
$m_\rho\
$exchange roles. This does not affect our results.
\endrem

\subhead {1.1} \endsubhead
For abelian varieties of types I, II and III  $\g_\ell,\ \h_\ell$ 
are semi-simple and have the same invariants on $V_\ell \otimes_\Ql V_\ell$ 
(cf.\ 0.4.3 and Remark 0.7.1(2)). In the case
of type IV, the Lie algebras can be non-semi-simple. Since $\g_\ell^{ss}
\subseteq \h_\ell^{ss},\ ${\it
a\ priori}\ $\g_\ell^{ss}$ can have more invariants than $\h_\ell^{ss}$ in
$V_\ell^{\otimes 2}$. However, in our special case the following is true.

\thm 
If $D = \k$, then $\g^{ss}_\ell$ and $\h^{ss}_\ell$ are non-self-dual.
\endst

\rem If $\g_\ell = \g_\ell^{ss},\ \h_\ell = \h_\ell^{ss}$ we get
nothing new (cf.\ 0.7.1).
\endrem

\pf  {\bf(0)}.  If $\overline \g^{ss}_\ell\ $is symplectic or
orthogonal, then so is $\overline
\g^{ss}_\ell \otimes_\Qlbar \Cl.\ $If $\overline \h^{ss}_\ell \otimes_\Qlbar
\Cl$ fixes a bilinear form
on $\overline V_\ell\ $coming from $V_\ell \otimes_\Ql \Qlbar$, then $\overline
\h^{ss}_\ell\
$fixes the form. So, we can extend scalars to $\Cl\ $and will use the same
notation $\overline
\h_\ell,\ \overline \g_\ell\ $for the corresponding extentions of the Lie
algebras.

{\bf(i)}.  Let us consider first the symplectic case.

If $\overline \g^{ss}_\ell$ {\it fixes} 1-dimensional subspace
$\chi \subseteq \overset 2 \to \wedge \overline V_\lambda$, then $\chi$ is a
1-dimensional
$\overline \g_\ell$-sub-\rep. The Hodge-Tate  decomposition implies
$$\overset 2 \to \wedge \overline V_\lambda = \overset 2 \to \wedge \overline
V_\lambda(0)
\oplus (\overline V_\lambda(1) \otimes \overline V_\lambda(0)) \oplus \overset
2 \to \wedge
\overline V_\lambda(1) .$$
The Hodge-Tate weight  of the terms on the right is 0, 1 and 2 respectively.

Since $\dim_\Cl(\chi) = 1,\ $it is of Hodge-Tate
weight 0, 1 or 2.
$\chi^{g/2} = \text{det}(\overline V_\lambda) := \overset g \to \wedge
\overline V_\lambda =
\overset {m_\sigma} \to \wedge \overline V_\lambda(1) \otimes \overset {m_\rho}
\to \wedge
\overline V_\lambda(0),\ $the Hodge-Tate weight  of the RHS is $m_\sigma.\
$Hence the Hodge-
Tate weight  of $\chi^{g/2}\ $is $m_\sigma.\ $On the other hand, the weight of
$\chi^{g/2}\ $is
$g/2$-times the weight of $\chi,\ $hence is equal to 0, $g/2\ $or $g.\ $So
$m_\sigma = 0, g/2\ $or
$g.\ $
The cases $m_\sigma = 0\ $or $g\ $correspond to the $\k$-signature $(0, g)\ $or
$(g, 0).\ $In
either case the \av\ is isogenous to a product of CM elliptic curves ([Sh~1,
Proposition
14]). In the case $m_\sigma = g/2$ we have $m_\sigma = m_\rho$ and the Lie
algebra
$\overline \g_\ell = \overline \g^{ss}_\ell$ (as well as $\overline \h_\ell =
\overline \h^{ss}_\ell)$
is non-self-dual (cf.\ 0.7.1), hence this $\chi\ $does not exist!

{\bf(ii)}.  The orthogonal case is a direct consequence of the symplectic and
0.4.3 (take $s = 2:\
\text{Sym}^2(\overline V_\lambda) \embed \overset 2 \to \wedge (2\overline
V_\lambda)).$ \QED
\endpf

\rem Another way to conclude that the \HT\ weight of $\chi$ is 1 is
to use a result of Raynaud that the Galois action on 1-dimensional
subrepresentations of (co)homology is by (powers of) the cyclotomic character
$\chi_\ell$. Either way, the result is a consequence of the
existence of the Hodge-Tate decomposition.
\endrem

\subhead {1.2} \endsubhead
We can apply this consideration of the Hodge-Tate
decomposition to abelian varieties
with $D = \k\ $and $\k$-signature $(m_\sigma, m_\rho)\ $such that
$\gcd(m_\sigma, m_\rho) = 1\
$(we call them  {\it Ribet-type}\ abelian varieties, cf.\ [Ri]), then an
argument of Serre
([S~1, $\S4$]) implies the Tate conjecture in that case. Indeed, Ribet's
proof of \loc.cit., Theorem 3, {\it verbatim}\ provides the Tate cycles are
generated by divisors and hence the following theorem.

\thm If $A$ is a Ribet-type abelian variety, then the Tate cycles (on
the abelian variety,
and all its self-products) are generated by divisors and hence the Tate, the
Hodge and the Mumford--Tate
conjectures hold. \QED
\endst

\rem Applicability of this argument to the Tate cycles was undoubtly known 
to Ribet and was also noticed in [C 1] and [MZ].
\endrem

\head
{\ff 2. Abelian 4-folds}
\endhead

\subhead {2.0} \endsubhead
The first non-trivial case of the Hodge, the Tate and the Mumford--Tate 
conjectures is that of abelian {\it fourfolds}. In that case 
$1 \leq (D:\Q) \leq 8$ and the dimensions of the irreducible sub\reps\ 
of $\h,\ \g_\ell$ are $\leq 8$.

The case of a simple 4-fold with $(D:\Q) \geq 2$ was studied in \cite{MZ}.
In this section we survey the 4-dimensional case just
indicating the ideas involved. For details (in particular, proof 
of Theorem 2.4 below) see \cite{MZ}.

\subhead {2.1} \endsubhead
If $(D:\Q) \geq 2$, then the dimensions $\leq 4$ and the restrictions 
imposed on $\g_\ell,\ \h$ (cf.\ 0.1.2, 0.3.1, 0.5) force the \reps\ to 
be {\it unique}, hence \avs\ in these cases verify MT. 

Moreover, in``most cases" these Lie 
algebras are the ``largest possible," viz., coincide with the Lie algebra 
of the Lefschetz group. Thus for any such \av\ (excluding type III, cf.\ 0.1.5) 
the Hodge and the Tate conjectures hold.

\subhead {2.2} \endsubhead
The case $D = \Q$ is slightly more subtle: the choice for 
$\g_\ell,\ \h$ is {\it not} unique anymore. In fact, both possibilities,
$\sl_2 \times \sl_2 \times \sl_2$ and $\sp_8$, do occur.

The first one, viz., $\overline \g_\ell = \overline \h_\ell \iso \sp_8$, 
is the generic case (cf.\ [Ab~1], [Ma]). Abelian varieties with 
$\overline \h \iso \sl_2 \times \sl_2 \times \sl_2$, were constructed by 
Mumford in [M 2]. 

However, imposing some {\it extra conditions} on a simple 4-fold with 
$D = \Q$ we still can conclude MT (cf.\ 5.2).

\subhead {2.2.1} \endsubhead
In the generic case 
$\overline \g_\ell = \overline \h_\ell = \overline \l_\ell$, where
$\overline \l_\ell := \Lie(L(A)) \otimes_\Q \Qlbar$, the Lie algebra of 
the Lefschetz group, and all the conjectures follow from 0.1.5.

\subhead {2.2.2} \endsubhead
In the Mumford case, $\overline \h \iso \sl_2 \times \sl_2 \times \sl_2$, 
again $\overline \g_\ell = \overline \h_\ell$ and MT holds. However, 
$\overline \l \iso \sp_8$ and {\it not} all the Hodge/Tate cycles (on 
{\it self-products}!) are divisorial.

\subhead {2.2.3} \endsubhead
Meanwhile, whatever the case,
$\sl_2 \times \sl_2 \times \sl_2$ or $\sp_8$,
[Ta, Lemma 4.10] implies that all the Hodge and the Tate cycles on the 
\av\ {\it itself} are divisorial. However, calculations show (cf.\ [H 1,
Lemma 5.2, (5.2.2)]) that on the ``square" of the Mumford 4-fold not all 
the Hodge cycles are divisorial.

\subhead {2.3.1} \endsubhead
If $D = \Q$, then either 
$\overline \g_\ell = \overline \h_\ell$
and MT holds, or 
$\overline \g_\ell \neq \overline \h_\ell$
and then 
$\overline \g_\ell \iso \sl_2 \times \sl_2 \times \sl_2,
\ \overline \h_\ell \iso \sp_8$.
In the latter case the (Lie algebra of the) Hodge group is equal
to the (Lie algebra of the) Lefschetz group, thus all the Hodge cycles 
on the self-products of the abelian variety are divisorial, see also $\S7$.

\subhead {2.3.2} \endsubhead
For the Weil case (cf.\ 0.7.3) {\it generically} the ring of Hodge
cycles is not generated by divisors ([W], also [MZ]). So, if there are any
doubts about the Hodge
conjecture (and hence the Tate conjecture), the Weil abelian varieties are the
ones to look at.
Recently, C. Schoen ([Sc], also [vG]) succeded in proving the Hodge
conjecture for one family of Weil
4-folds admitting an action of $\Q(\mu_3).$

\subhead {2.3.3} \endsubhead
2.1 answers a question of Tate (cf.\ [T~2], p.~82) on whether
the Tate conjecture is
true for the Schoen family. 

\subhead {2.4} \endsubhead
We summarize the above discussion in the following theorem.

\thm 1.  If $A$ is any 4-dimensional abelian variety, then the
rings of the Tate cycles and
the Hodge cycles coincide (hence, the Hodge and the Tate conjectures for this
variety are equivalent).

2.  If, additionally, $\End^{\circ}(A) \neq \Q,\ $then MT holds.
\endst

\rems 1.  Later (Theorem 5.2) we will see that even when
$\End^{\circ}(A) = \Q$, MT holds under some reduction conditions.

2.  Recall (0.4) that MT implies that the Hodge and the Tate conjectures are
equivalent for an \av\ and {\it all its self-products.}

3.  If $A$ is non-simple, then the second hypothesis of the theorem
is satisfied.
\endrem

\pf The only ``real" case to consider is that of  $A$ {\it simple}. 
The proof can be found in [MZ] (in fact, Moonen-Zarhin
considered a deeper problem of ``when and why" a simple 4-fold has an
exceptional Weil class).

If $A$ is a {\it non}-simple abelian 4-fold, say $A$ is isogenous to $A_1
\times A_2$, then
$\dim(A_i) \leq 3$. Hence, by the (1,1)-theorems and duality,
all the Hodge cycles
on $A_i$ are divisorial. The embeddings $\g_\ell\subset \h_\ell
\embed \sp(V_\ell)$
factor through the sub-\reps\ corresponding to the simple components of $A$.
The dimensions of the sub-\reps\ are $\leq 6$, and there is not ``enough room" 
for $\g_\ell$ and $\h_\ell$ to be different, i.e.,
$\g_\ell = \h_\ell$, hence MT holds. \QED
\endpf

\subhead {2.5} \endsubhead
Let us indicate what is the situation regarding the Hodge and the
Tate conjectures for
{\it non-simple}\ abelian 4-folds. As above, let $A$ be isogenous to
$A_1 \times A_2$. Then all the Hodge and the Tate cycles on the $A_i$'s are
divisorial.

\subhead {2.5.1} \endsubhead
We can also say that all the Hodge cycles (and hence the Tate
cycles) on $A$ and all its self-products are
generated by divisors in the following cases:
\roster
\item"{\bf 1. }"  Neither of the $A_i$'s is of type IV ([H~2, Theorem 0.1]).
\item"{\bf 2. }"  $A_1\ $is not of type IV, $A_2\ $is of CM-type (\loc.cit.,
Proposition 3.1).
\item"{\bf 3. }"  If the $A_i$'s are {\it non}-CM, type IV abelian surfaces, then
according to [Sh~1, Theorem 5, Propositions 17, 19], the $A_i$'s are products
of CM elliptic curves. Hence so is
$A = A_1 \times A_2\ $and for such \avs\ the result stated above is known
([Im]; [H~1,
Theorem 2.7]).
\item"{\bf 4. }"  If the $A_i$'s are {\it isogenous}\ CM surfaces, then 
by remark 0.1.2(1) and 0.1.5, $\Hg(A) = \Hg(A_1),\ L(A) = L(A_1)$. 
By 0.1.5(1) $L(A_1) = \Hg(A_1)$, hence $L(A) = \Hg(A)$, and applying 
0.1.5(1) once again we conclude the result.
\endroster

\subhead {2.5.2} \endsubhead
For the remaining case, viz., both the $A_i$'s are {\it non-isogenous}
CM abelian varieties, let us just mention that Shioda constructed an 
example of a product of a simple CM 3-fold,
$A_1$, with a CM elliptic curve, $A_2$, \s.t. on $A = A_1 \times A_2$
there exist
exceptional, non-divisorial, Hodge cycles, [Shi, Example 6.1]. In this example,
however, the
Hodge (hence the Tate) conjecture holds.

\head
{\ff II. Abelian varieties with reduction conditions}
\endhead
\head
{\ff 3. Bad reduction and monodromy action}
\endhead

\subhead {3.0} \endsubhead
Let $A\ $be an abelian variety defined over a number field $K$.
Assume $A\ $has bad reduction at a prime $\wp\ $of $\O_K.$ Let
$\tilde A\ $be the identity component of the special fiber of the N\'eron model
of $A.\ $Then $\tilde A\ $is semi-abelian:
$$0 \to H \to \tilde A \to B \to 0,$$
where $H\ $is the affine subgroup of $\tilde A,\ B\ $is the abelian
quotient.

\subhead {3.0.1} \endsubhead
Since we are concerned with the Lie algebra of the (image of)
Galois, we may pass
to a finite extension of $K\ $(cf.\ Remark 1 in 0.3). So, according to the
semi-stable reduction theorem ([G, Th\'eor\`em 3.6]), by extending the base
field if necessary,
we may assume that the reduction is {\it semi-stable}\ (i.e., $H$ is a torus)
and {\it split}\ (i.e., $H\ $is split: $H \iso {\Bbb G^r_m})$.

The dimension $r$ of $H$ we call the {\it toric\ rank} of (the reduction of)
$A$.

\subhead {3.0.2} \endsubhead
$D = \End^{\circ}(A)$ as before, there is a homomorphism $D \to
\End^{\circ}(H)$, 
$1_A \mapsto 1_H$. But $\End^{\circ}({\Bbb G^r_m}) = M_r(\Z)
\otimes_\Z \Q$, hence $(D:\Q) | r.$

\subhead {3.1} \endsubhead
Consider the corresponding ``specialization sequence"
$$
0 \to V_\ell^\I \to V_\ell(A) \to U \to 0,
$$
where $ V_\ell(A)\ $is the Tate module of $A$,
$\I := \I(\wp)$ is the inertia group at $\p$,
$V_\ell^\I := V_\ell(A)^\I$ is the submo\-dule of $\I$-invariants and $U$ is a
trivial $\I$-module (cf.\ [G, Proposition 3.5]).
We have $\dim_\Ql(V_\ell(A)) = 2g,\
\dim_\Ql(V_\ell^\I) = 2g - r,\ \dim_\Ql(U) = r$.

\subhead {3.2} \endsubhead
The above sequence is a sequence of of $\I$-modules.
The $\I$-action is called the {\it local\ monodromy}\ action.

\subhead {3.3} \endsubhead
The reduction map at $\wp$ induces an isomorphism $V^\I_\ell \iso
V_\ell(\tilde A)$, the Tate module of $\tilde A$, and takes a submodule $W\!
\subseteq\! V^\I_\ell$ to the Tate module $V_\ell(H) \subseteq V_\ell(\tilde
A)$ of the toric part $H$ of $\tilde A$ (cf.\ [ST, Lemma 2], [G, 2.3], [I],
[O]). In fact, according to the ``Igusa-Grothendieck Orthogonality Theorem," $W =
(V^\I)^{\perp}$ \wrt\ the Weil pairing on $V$ (cf.\ [I, Theorem 1], 
[G, Th\'eor\`eme 2.4], also [O, Theorem (3.1)]).

\subhead {3.4} \endsubhead
The monodromy action on $V_\ell(A)\ $is, in general,
quasi-unipotent (e.g., [G], [ST], [O]). However, since (we assumed that) the
reduction of $A$is {\it semi-stable}\ and {\it split},\ this action is, 
in fact, {\it unipotent} (cf.\ [G, Corollaire 3.8]).

\subhead {3.4.0} \endsubhead
Picking a vector subspace $T$ of $V_\ell(A)\ $specializing to
$U$, we get the matrix form of the monodromy action:
$$\aligned V^\I_\ell &\left \{  {\aligned &W\{\\ &\phantom{T\{} \endaligned}
\right.\\
		      T &\{ \endaligned
\pmatrix 1_r & 0 & *_r\\
\\
0 & 1_{2g-2r} & 0\\
\\
0 & 0 & 1_r
\endpmatrix.
$$
\subhead {3.4.1} \endsubhead
Passing to the Lie algebra $\i := {\Lie}(\I),\ $we conclude the
existence of nilpotents, $\tau
\in \i \subset \g_\ell,\ $of order 2, i.e., $\tau^2=0$, and rank (\wrt\ $V_\ell)$
$ \rk_{V_\ell}(\tau) \leq r$,\ where $\rk_{V_\ell}(\tau) := 
\dim_\Ql(\tau V_\ell)$\ = rank of the matrix of $\tau \in \gl(V_\ell)$. 

\subhead {3.4.2} \endsubhead
The Neron-Ogg-Shafarevich criterion ensures that $\exists\, \tau \neq 0$,
since $A$\ has bad reduction.

\subhead {3.4.3} \endsubhead
Moreover, if $N$\ is given by the above matrix, then $\tau = N - 1_{2g}$ 
is the logarithm of the
monodromy corresponding to the {\it monodromy\ filtration}\ (cf.\ [G,~4.1 and
also Corollaire 4.4]; also [Il, 2.6])
$$
(0) \subset W \subset V^\I_\ell \subset V_\ell,
$$
and $\tau$ maps $V_\ell \to W,\ V^\I_\ell \to 0$, inducing an isomorphism 
of the quotients $V_\ell/V_\ell^\I \overset{\sim} \to \to W$, or, 
in our notation,
$$
\tau: T \overset{\sim} \to \to W
$$
(cf.\ [G, 4.1.2], [Il, (2.6.3)]; see also our Remark in 0.3 for why we omit
the Tate twist in this formula). In particular, $\rk_{V_\ell}(\tau) = r$.

\subhead {3.5} \endsubhead
By extending scalars to $\Qlbar\ $we get the corresponding
nilpotents (of the same
order) in each \irr\ component of $V_\ell \otimes \Qlbar$ with the sum of
the ranks \wrt\ each of the components being equal to the rank \wrt\
$V_\ell$.

\head
{\ff 4. Minimal reduction}
\endhead

\subhead {4.0} \endsubhead
We say that an abelian variety $A\ $over a \nf\ has {\it minimal}\ bad
reduction at a prime
$\wp$ of this field (or, just minimal reduction, for short) if the reduction
is bad and the rank of the
toric part $H\ $of $\tilde A\ $(cf.\ 3.0) is the minimal possible.

\subhead {4.1} \endsubhead
Let us go back to the case $D = \k,\ $in which $\overline
\g_\ell \subseteq \overline \h_\ell
\subseteq \gl(W_\lambda) \overset \Delta \to \to$
$ \sp(W_\lambda
\oplus \check W_\lambda)\ $(cf.\
0.6.2). ~The toric rank should be even (cf.\ 3.0.2), say, $2r.\ $If
$\tau'=\Delta(\tau) \in \Delta(\overline
\g_\ell)\ $is a nilpotent of rank $\rk_{V_\ell \otimes \Qlbar}(\tau') =
2r$ (cf.\ 3.4.3), then $\tau^2 = 0,\ \rk_{W_\lambda}(\tau) = r.$

\subhead {4.1.1} \endsubhead
3.0.2 and 3.4.3 imply that in the case $D = \k$ the minimal toric rank is 2.

\subhead {4.1.2} \endsubhead
Hence in the minimal reduction case $\exists \, \tau \in
\overline \g_\ell \subseteq
\overline \h_\ell \subseteq \gl(W_\lambda)$ such that $\tau^2 = 0,\
\rk_{W_\lambda}(\tau) = 1$.\ The same, clearly, holds if we replace 
the Lie algebras with their \ss\ components, since
all nilpotents live in these components.

Let 
$U := W_\lambda$, and (as in 0.4.4)
$\gg = \overline \g_\ell^{ss},\ \hh = \overline \h_\ell^{ss}$,
and rewrite the above as:
$$
\gather
\tau \in \gg \subseteq \hh \subseteq \sl(U),\ \tau^2 = 0,\ \rk_U(\tau) =
1,\\
\gg,\ \hh\ :\ \text{\ss\ \irreps}.
\endgather
$$

Such an element $\tau\ $of rank 1 and order 2 is called a {\it transvection}.

\subhead {4.2} \endsubhead
It is a very restrictive condition for an \irrep\ of a
semi-simple Lie algebra to contain a transvection.

\lemma If $\a \embed \sl(U)\ $is a \ss\ faithful \irrep,\ $\tau \in \a,\
\tau^2 = 0,\ \rk_U(\tau) = 1$, 
then $\a$ is simple and, moreover, it is either $\sp(U)\ $or $\sl(U).$
\endst

\pf This is proved in [McL] (cf.\ also [PS]).\QED
\endpf

\subhead {4.3} \endsubhead
We will also need the following simple fact.

\lemma If $\rk_U(\tau)$\ is prime to $\dim(U)$, then $\a$ is simple. \QED
\endst

\head
{\ff 5. Applications of minimal reduction}
\endhead

\subhead {5.1} \endsubhead
An immediate application of 4.1.2 (existence of rank 1 quadratic
nilpotents in $\g$) and Lemma 4.2 is the following theorem.

\thm If $A$ is a simple abelian variety with $D \subseteq \k$,
having minimal reduction, then MT holds.
Moreover, if $D = \Q,\ $then all the Hodge and the Tate cycles are divisorial,
hence the Hodge and the Tate conjectures hold.
\endst

\pf 1. If $D = \Q,\ $then the Tate module $V_\ell(A)\ $is absolutely
\irr\ and symplectic.
The minimality of reduction implies that the rank of a correspondent nilpotent
is 1. The result now follows from 4.2.

2.  If $D = \k,\ $the result follows from 4.1.2, 4.2 and 1.1 (cf.\ 0.1.2). \QED
\endpf

\rems 1.  Such abelian varieties exist and, moreover, form a
subset dense in the complex topology in the corresponding moduli space (cf.\
[L]).

2.  The importance of the Weil type abelian varieties is not limited to the
fact that they (may)
have non-divisorial Weil cycles (cf.\ 4.4). Proving the algebraicity of the Weil
cycles is a critical
ingredient in proving the Tate conjecture 
(cf.\ [D~3, $\S\S\ $4$\sim$6]; also [An], [Ab~2, $\S\ 6$]).
\endrem

\subhead {5.2} \endsubhead
Using the same method we can now extend Theorem 2.4 in the following way.

\thm  If $A\ $is a simple abelian 4-fold with $D = \Q$ admitting bad
but not purely
multiplicative reduction, then all the Hodge and the Tate cycles are
divisorial, hence the Hodge
conjecture, the Tate conjecture and MT hold.
\endst

\pf  The possible values of the toric rank in this case are 1, 2, 3 (4
corresponds to the
purely multiplicative reduction). 

First recall (2.2) that the only choices for $\gg$ and $\hh$ are $\sp_8$ or
$\sl_2 \times \sl_2 \times \sl_2$.
But $\sl_2 \times \sl_2 \times \sl_2$ does not contain quadratic nilpotents of
rank 1, 2 or 3.
So, $\gg = \hh \iso \sp_8$ and 0.1.5 finishes the proof.  \QED
\endpf

\rems 1.  Applying results/methods of [Z~5], [LZ] one can also verify all the 
conjectures for a simple abelian 4-fold admitting certain types of
{\it good} reductions, further restricting the class of 4-folds for which the
conjectures are not yet known.


2.  G. Mustafin [Mus] has proved a ``geometric analog" of MT for (families of)
abelian varieties with purely multiplicative
reduction (i.e., the algebraic envelope of the image of the monodromy coincides
with the Hodge group for a ``sufficiently general" \av\ in the family), cf.\
also [H~3]. It appears, though, that his methods cannot be transplanted to the
arithmetic situation.
\endrem

\head
{\ff 6. Another type of bad reduction}
\endhead

\subhead {6.0} \endsubhead
Now we want to establish a result analogous to Theorem 5.1 for
another type of bad reduction.
Namely, consider an \av\ $A$ admitting bad (semi-stable, split) 
reduction of toric rank $r$ \s.t.\ $r/(D:\Q)$ is prime to $2\dim(A)/(D:\Q)$, 
where as before $D = \End^\circ(A)$.
Assume also $D$ is {\it commutative}. In this case $\gg,\ \hh$\ (notation as in
4.1.2) are simple (cf.\ 4.3) and contain nilpotents of rank prime to the
dimension of the representations (cf.\ 3.4.3). In place of Lemma 4.2 we use the
results of Premet-Suprunenko
[PS] on classification of quadratic elements (= nilpotents of order 2) and
quadratic modules
(= representations containing non-trivial quadratic elements) of simple Lie
algebras.
\footnote"$^\dagger$"{This terminology is apparently standard in the finite
groups theory, cf.\ [Th].}

We use the fact that the representations of $\gg,\ \hh$ are minuscule and $\gg$ 
is not exceptional (0.5.1).

We may assume that the dimensions of the representations is $> 4.$

\subhead {6.1} \endsubhead
One of our key tools replacing Lemma 4.2 is the following result.
\thm If $\a \subset \b \embed \sl(U),\ \a \neq \b$, both Lie
algebras are simple and the
representations are (faithful) irreducible and minuscule, then $\b$ is
classical and (its highest weight is) $\varpi_1.$
\endst

\pf Since any minuscule \rep\ is quadratic (cf.\ [B, Ch VIII, $\S$7.3]),
we can apply [PS, Theorem 3], and exclude non-minuscule cases. \QED
\endpf

\subhead {6.2} \endsubhead
We will be interested in 2 cases: $D = \k$ and $D = \Q$. In the former case we
know that $\gg$ and $\hh$ are both non-self-dual (cf.\ Theorem 1.1) and if they
satisfy the conditions of the theorem, then
\roster
\item"{\bf --}"
$\hh$ is classical, $\varpi_1$ and non-self-dual, hence $\hh = (A_n, \varpi_1)$
\item"{\bf --}"
$\gg$ is classical, minuscule and non-self-dual, hence $\gg = (A_m, \varpi_s),\ m
\neq 2s$.
\endroster

If $D = \Q$, we know that $\gg,\ \hh$ are symplectic and we again
have a unique possibility for $\hh$, viz., $\hh = (C_n, \varpi_1)$. However,
there are several {\it a priori} possible choices for $\gg$.

To eliminate (as many as we can) possibilities of $\gg \neq \hh$ we use the
existence in the representations of a quadratic nilpotent of {\it rank prime to
the dimension} of the \rep.

So, we consider a slightly more general situation. As above
\roster
\item"$\bullet$"
$\gg \varsubsetneq \hh \subset \sl(U),\ \gg,\ \hh$ are classical simple Lie
algebras,
\item"$\bullet$"
$\gg$ is minuscule, $\hh$ is $\varpi_1$,
\item"$\bullet$"
there exists $\tau \in \gg \subset \hh \subset \sl(U)$ with $\tau^2(U) = 0$,
$\rk_U(\tau) = r,\ \dim(U) = n$ and $\gcd(r, n) = 1$.
\endroster

We add the following condition wich is satisfied in both our cases:
\roster
\item"$\bullet$"
$\gg,\ \hh$ are non-self-dual, or orthogonal, or symplectic {\it simultaneously}.
\endroster

\subhead {6.3} \endsubhead
First we exclude the cases $\gg = (D_m, \varpi_{m-1}),\ (D_m,
\varpi_m)\ $for $m > 4.$

\lemma If $\tau \in \gg \embed \sl(U),\ \gg = (D_m, \varpi_{m-1})$ or
$(D_m, \varpi_m),\ \tau$ is a quadratic element, then $\gcd(r,n) > 2.$
\endst

\pf [PS, Lemma 21, Note 2, Lemma 17] imply $r = 2^{m-3}\ $or $2^{m-2}$,
while $n = 2^{m-1}$. \QED
\endpf

$\gg \neq (D_4, \varpi_3),\ (D_4, \varpi_4)$ either. It is enough
to show this for $\varpi_4$, since they are (graph)isomorphic.

\prop $(D_4, \varpi_4) \not\embed\ $(classical, $\varpi_1).$
\endst

\pf Note that $n$ = 8 in here. We do this case by case:
\roster
\item"\bf 1.  "
${(D_4, \varpi_4) \not\embed\ (B_\bullet, \varpi_1)}$, since
the dimension of the
RHS is odd.

\item"\bf 2.  "
${(D_4, \varpi_4) \not\embed\ (C_\bullet, \varpi_1)}$, since the
LHS is orthogonal
while the RHS is symplectic (cf.\ [B, table 1]).

\item"\bf 3.  "
${(D_4, \varpi_4) \not\embed\ (D_4, \varpi_1)}$ (e.g.,
[Z~2, $\S5,\ $Key lemma], although this is overkill).

\item"\bf 4.  "
$(D_4, \varpi_4)\not\embed (A_\bullet, \varpi_1)$, since the 
LHS is orthogonal, the RHS is not. \QED
\endroster
\endpf

\prop If $\gg = (D_m, \varpi_1)$, then $\gg = \hh.$
\endst

\pf $\hh = A_\bullet,\ C_\bullet$ are excluded: $\gg$ is orthogonal, $\hh$
is not; $\hh = B_\bullet$ is excluded by a dimensional reason. \QED
\endpf

\prop If $\gg = (C_m, \varpi_1)$, then $\gg = \hh.$
\endst

\pf $\hh = A_\bullet,\ B_\bullet,\ D_\bullet$ are not symplectic ... \QED
\endpf

\prop If $\gg = (B_m, \varpi_1)$, then $\gg = \hh.$
\endst

\pf $\hh = A_\bullet,\ C_\bullet$ are not orthogonal...
$(D_\bullet, \varpi_1)$ is even-dimensional...  \QED
\endpf

\prop If $\gg = (A_m, \varpi_s)$, then $\hh$ must be $(A_{n-1}, \varpi_1).$
\endst

\pf First note that the only self-dual representation of $A_m$ is
$\varpi_s$ with $s = {{m+1} \over 2}$ (there is no such a representation 
if $m$ is even). So, if $\hh$ is self-dual, then
so is $\gg$ (cf.\ 6.2) and we may assume $s = {{m+1} \over 2},\ m$ : odd. But
then $\dim(A_m, \varpi_s) = \binom{2s}s$ : even, thus $\hh \neq B_\bullet$.
To exclude the other cases (i.e., $C_\bullet,\ D_\bullet$) we use the fact 
that $r = \rk(\tau) = \binom{m-1}{s-1}$ ([PS, $\S2\ \&$ Lemma 18]):
$$
r = \binom{2(s-1)}{s-1},\quad
n = \binom{2s}{s} = \binom{2(s-1)}{s-1}{{(2s-1)2s} \over {s^2}} = r{{2(2s-1)}
\over s};
$$
$\gcd(n,r) = 1 \imply r|s$, which is not true: $n > 4 \imply s > 3 \imply
\binom{2(s-1)}{(s-1)} > s$.
\QED
\endpf
%
%
%

\rem If the $\k$-signature of the abelian variety is $(m_\sigma, m_\rho)\
$with $m_\sigma
\neq m_\rho,\ $then, even if $\gg$ is not simple, simple components of
$\gg$ are of type $A$ (cf.\ [Y]).
\endrem

So, the only possibility for $\gg \subsetneqq \hh$ is $\gg = (A_m, \varpi_s)$
for some $s,\ \hh = (A_{n-1}, \varpi_1)$. In this case we can say the following.

\prop Let $\gg = (A_m, \varpi_s) \embed \hh = (A_{n-1}, \varpi_1) \iso \sl(U)$
(fix the isomorphism), $2 \leq s < {{m+1} \over 2},\ \tau \in \gg,\ r =
\rk_U(\tau),\ \gcd(n,r) = 1$. Then either $s = 3,\ m = 7$, or $s = 2$.
\endst

\pf $r = \binom{m-1}{s-1}\ $([PS, $\S2\ \&\ $Lemma 18]), $n = \binom{m+1}
s = \binom{m-1}{s-1}{{m(m+1)} \over {s(m+1-s)}}, \gcd(r,n) = 1 
\imply r\ |\ s(m+1-s)$
and the result follows from the following simple observation:

$\binom{m-1}{s-1} | s(m+1-s)$ if and only if $(m, s) = (7,3)$ or $s = 2$.\QED
\endpf
%
%
%
%

\remark
{Remark}
If $s = 2$, then $r = \binom{m-1}{s-1} = m - 1,\ 
n = \binom{m+1} s = {{m(m+1)} \over 2}$. 
Since $\gcd(m, m-1) = 1$ and $\gcd(m-1, m+1)$ = 1 or 2, $\gcd(r, n) = 1$ 
if and only if $m$ is even or $m \equiv 1\ \pmod 4$ 
(i.e., $m \not\equiv 3\ \pmod 4)$. 
\endremark

\subhead {6.4} \endsubhead
So, if $D = \k$ we have the following result.

\thm If $A$ is a simple abelian variety with $D = \k,\ g = \dim(A)$ having bad
reduction with the toric rank $2r$ and $r$ is prime to $g$, then MT holds if
$(g, r)$ is neither
(56, 15) nor of the form $({{m(m+1)} \over 2}, m-1)$.
\QED
\endst

\subhead {6.5} \endsubhead
If $D = \Q$, then $\gg,\ \hh$ are symplectic, hence the
theorem holds with no exceptions.

\thm If $A$ is an abelian variety with $D = \Q,\ g=\dim(A)$ having
bad reduction with the toric rank $r$ prime to $2g$, then MT holds. \QED
\endst

\subhead {6.6} \endsubhead
Using the same methods one can handle the case of quadratic
elements of rank 2. Namely, the following result holds.

\thm Let $A$ be a simple abelian variety with $D = \Q$. If $A$ has bad
reduction with toric rank 2, then MT holds.
\endst

\pf (Sketch) Let $\dim(A) = g$. We can assume $g \ge 4$.

One can check (cf.\ [PS, $\S2$, Lemma 18]) that if $\b \embed \sl(U)$ is a \ss\
irreducible classical Lie algebra,
$\dim(U) \ge 8$ and $\exists \, \tau \in \b,\ \tau^2 = 0,\ \rk_U(\tau) = 2$,
then either $\b$ is simple, and hence $\b = \sl(U),\ \sp(U)$ or
 $\so(U)$ (since $\dim(U) \ge 8,\ \b \not\iso (A_3, \varpi_2)$),
or $\b = \a \times \sl_2$, where $\a \iso \sl_g$ or $\sp_g$.
Since $\gg,\ \hh$ are {\it symplectic}, the only possibility
is $\sp(U)$ ! \QED
\endpf

\subhead {6.7} \endsubhead
All the varieties considered in $\S 6\ $exist and dense in the
(complex topology in the) corresponding moduli spaces (cf.\ [L]).

\subhead {6.8} \endsubhead
1.  The idea of using special element(s) in the representation
of the Hodge
group has been used before. However, to our knowledge, in those earlier cases
the element was
semi-simple of low rank (e.g., [Z~1]) and the results then follow from a
theorem of
Kostant [Ko] (cf.\ also [Z~3]).

2.  Katz used special unipotent elements to show that certain monodromy groups
are large.
However, the unipotents he considered were of the maximal possible rank, i.e.,
having only {\it one}\ Jordan block ([Ka~1, Ch. 7]).
\footnote"$^\ast$"{As an application of Katz's classification of
representations containing such unipotents (cf.\ \loc.cit., 11.5$\sim$11.7) 
one can
find modular curves for which the image of Galois in the corresponding
$\ell$-adic \rep\ is large.}

Katz was also using \ss\ elements for similar purposes, [Ka~2].

\head
{\ff 7. A curious result}
\endhead

\subhead {7.1} \endsubhead
Let us mention another application of Theorem 6.1.
\thm If $A$ is a simple abelian variety with (semi-simple parts of) 
$\h,\ \g_\ell$ {\it simple}, satisfying one of the following conditions:
\roster
\item"-"  the variety is of type I or II,
\item"-"  $D = \k$, the variety is of non-Weil type (i.e., the $\k$-signature is
$(m_\sigma, m_\rho)$ with $m_\sigma \neq m_\rho$), 
\endroster
then one of the following must hold:
\roster
\item"-" MT holds for $A$, 
\item"-" all the Hodge cycles (on $A$ and all its self-products) 
are generated by divisor classes (hence the Hodge conjecture holds). 
\endroster
\endst

\pf  As we mentioned in 6.0, the representations of $\h_\ell,\ \g_\ell$
are minuscule,
hence quadratic (cf.\ [B, Ch VIII, $\S$7.3, Proposition 7]), then so is 
$\l_\ell := \l \otimes_\Q \Ql$, where $\l$ is the Lie algebra of the 
Lefschetz group (cf.\ 0.1.5), and
$\gg \subset \hh \subset \ll \subseteq \sl(U)$ (cf.\ 0.2.3; here 
$\ll := \overline \l^{ss}_\ell,\  \overline \l_\ell 
= \l_\ell \otimes_\Ql \Qlbar$). If $\gg \subsetneqq \hh$ (i.e., 
MT does {\it not} hold), then by Theorem 6.1 $\hh \embed \sl(U)$ is
classical and $\varpi_1$. Thus $\ll$ is also simple, classical and $\varpi_1$.
We want to show that in this case $\overline \h_\ell = \overline \l_\ell$ and
the theorem then follows from 0.1.5.

Consider first the case of an abelian variety of type I or II. The Lie
algebras $\overline \h_\ell$ and $\overline \l_\ell$ are both symplectic,
simple, classical and $\varpi_1$. Hence $\overline \h_\ell 
= \sp(U)\ (\text{resp.}\ \so(U)) =  \overline \l_\ell$.

If $D = \k$, then $C_{\overline \h_\ell}$\ is 1-dimensional (0.7.2),
hence $C_{\overline \h_\ell} = C_{\overline \l_\ell}$ and 
$\hh = \ll = \sl(U)$ (cf.\ 1.1; also [Mu, Lemma 2.3]). The theorem follows. \QED
\endpf

\subhead {7.2.  Remarks} \endsubhead
1.  As one can see from the proof, if there is a way to
assure that $\h$ is simple and $\varpi_1$, then $\h=\l$, hence 
(for types I or II) the Hodge conjecture holds (e.g., 2.3).

2.  Abelian varieties with bad reduction as in 6.4$\sim$6.5 have simple
(semi-simple parts of) Hodge and Galois groups (cf.\ 4.3).

3.  For the Weil type varieties, i.e., $m_\sigma = m_\rho,\  C_{\overline
\h_\ell} = \{0\} \neq C_{\overline \l_\ell}$\ {\it generically}.\ 
This is (a restatement of) the main result of [W].

\enddocument